\documentclass[12pt]{article}

\usepackage[english]{babel}
\usepackage{upref,amsfonts,amsxtra,a4wide,latexsym}
\usepackage{amssymb,a4wide,doublespace,amsmath,epsf,mathrsfs,theorem}
\usepackage[matrix,arrow,cmtip,rotate]{xy}
\newdir{ (}{{}*!/-5pt/\dir^{(}}
\newdir{ )}{{}*!/-5pt/\dir_{(}}

\newtheorem{theorem}{Theorem}[subsection]
\newtheorem{lemma}[theorem]{Lemma}
\newtheorem{corollary}[theorem]{Corollary}

\newtheorem{proposition}[theorem]{Proposition}

{\theorembodyfont{\upshape}

}

\numberwithin{equation}{section}
\numberwithin{theorem}{section}

\newcommand{\Ker}{\mathrm{Ker}}
\newcommand{\image}{\mathrm{Im}}
\newcommand{\Dom}{\mathrm{Dom}}

\newcommand{\diag}{\mathrm{diag}}
\newcommand{\Id}{\mathrm{id}}

\newcommand{\Span}{\mathrm{span}}

\newcommand{\Tr}{\mathrm{Tr}}

\newcommand{\setmin}{\! \setminus \!}
\newcommand{\Ca}{\Omega}

\newcommand{\cM}{{\cal M}}
\newcommand{\ep}{{\varepsilon}}

\newcommand{\cC}{{\cal C}}
\newcommand{\cD}{{\cal D}}

\newcommand{\cB}{{\cal B}}

\newcommand{\cF}{{\cal F}}
\newcommand{\cO}{{\cal O}}

\newcommand{\A}{{\cal A}}
\newcommand{\C}{{\mathbb C}}
\newcommand{\Z}{{\mathbb Z}}
\newcommand{\N}{{\mathbb N}}

\newcommand{\I}{{\ensuremath{[0,1)}}}
\newcommand{\J}{{\ensuremath{[0,1]}}}

\newcommand{\cK}{{\cal K}}

\newcommand{\R}{{\mathbb R}}
\newcommand{\Cs}{{$C^*$-al\-ge\-bra}}

\newcommand{\gange}{\! \cdot \!}

\newcommand{\sh}{{$^*$-ho\-mo\-mor\-phism}}

\newcommand\eqdef{{\;\;\overset{\mbox{\scriptsize def}}{=}\;\;}}

\newenvironment{proof}[1][Proof:]%
{\begin{trivlist}\item[]\textbf{#1} }%
{\hbox{}\nobreak\hfill\quad\hbox{$\square$}\end{trivlist}}

\begin{document}
 
\title{A purely infinite AH-algebra and an application to AF-embeddability}
\author{Mikael R\o rdam}
\date{}
\maketitle

\begin{abstract} \noindent We show that there exists a purely infinite
  AH-algebra. The AH-algebra arises as an inductive limit of \Cs s of
  the form $C_0(\I,M_k)$ and it absorbs the Cuntz algebra $\cO_\infty$
  tensorially. Thus one can reach an $\cO_\infty$-absorbing \Cs{} as an
  inductive limit of the finite and elementary \Cs s
  $C_0(\I,M_k)$. 

As an application we give a new proof of a recent theorem of Ozawa
that the cone over any separable exact \Cs{} is AF-embeddable, and we
exhibit a concrete AF-algebra into which this class of \Cs s can be
embedded. 
\end{abstract}

\section{Introduction}
\noindent Simple \Cs s are divided into two disjoint subclasses: those
that are stably finite and those that are stably infinite. (A
simple \Cs{} $A$ is stably 
infinite if $A \otimes \cK$ contains an infinite projection, and it is
stably finite otherwise.) All simple, stably finite \Cs s admit a
non-zero quasi-trace, and all exact, simple, stably finite \Cs s admit a
non-zero trace.

A (possibly non-simple) \Cs{} $A$ is in \cite{KirRor:pi} defined to be
\emph{purely infinite} if no non-zero quotient of $A$ is abelian and
if for all positive elements $a,b$ in $A$, such that $b$ belongs to the
closed two-sided ideal generated by $a$, there is a sequence $\{x_n\}$
of elements in $A$ with $x_n^*ax_n \to b$. Non-simple purely infinite
\Cs s have been investigated in \cite{KirRor:pi}, \cite{KirRor:pi2}, and
\cite{BlanKir:pi3}. All simple purely infinite
\Cs s are stably infinite, but the opposite does not hold, cf.\
\cite{Ror:simple}. 

The condition on a (non-simple) \Cs{} $A$, that all projections in $A
\otimes \cK$ 
are finite, does not ensure existence of (partially defined)
quasi-traces. There are stably projectionless purely infinite \Cs s---
take for example $C_0(\R) \otimes \cO_\infty$, where 
$\cO_\infty$ is the Cuntz algebra generated by a sequence of
isometries with pairwise orthogonal range projections---and purely
infinite \Cs s are traceless. 

That stably projectionless purely infinite \Cs s can share
properties that one would expect 
are enjoyed only by finite \Cs s was demonstrated in a recent paper
by Ozawa, \cite{Oza:AF}, in which is it shown that the suspension and
the cone over any separable, exact \Cs{} can be embedded into an
AF-algebra. (It seems off hand reasonable to characterize
AF-embeddability as a finiteness property.) In particular,  $C_0(\R)
\otimes \cO_\infty$ is 
AF-embeddable and at the same time purely infinite and
traceless. It is surprising that one can embed a traceless
\Cs{} into an AF-algebra, because AF-algebras are well-supplied with
traces. If
$\varphi \colon C_0(\R) \otimes \cO_\infty \to A$ is an embedding into
an AF-algebra $A$, then $\image(\varphi) \cap \Dom(\tau)
\subseteq \Ker(\tau)$ for every trace $\tau$ on $A$. This can happen
only if the ideal lattice of $A$ has a sub-lattice isomorphic to the
interval $[0,1]$ (see
Proposition~\ref{prop:manyideals}). In particular $A$ cannot be simple. 

Voiculescu's theorem, that the cone and the suspension over any
separable \Cs{} is quasi-diagonal, \cite{Voi:qd}, is a crucial
ingredient in Ozawa's proof.

By a construction of Mortensen, \cite{Mor:IO2}, there is to
each totally ordered, compact, metrizable set $T$ an AH-algebra $\A_T$
with ideal lattice $T$ (cf.\ Section~2). A \Cs{} is an AH-algebra, in
the sense of Blackadar \cite{Bla:matricial_topology}, if
it is the inductive limit of a 
sequence of \Cs s each of which is a direct sum of \Cs s of the form
$M_n(C_0(\Omega)) = C_0(\Omega,M_n)$ (where $n$ and $\Omega$ are
allowed to vary). We show in Theorem~\ref{thm:A} (in combination with
Proposition~\ref{prop:A=B}) that the AH-algebra $\A_\J$ is
purely infinite (and hence traceless)---even in the strong sense that it
absorbs $\cO_\infty$, i.e., $\A_\J \cong \A_\J \otimes \cO_\infty$---
and $\A_\J$ is an inductive limit of \Cs s of 
the form $C_0(\I,M_{2^n})$. We can rephrase this result as follows:
Take the smallest class of \Cs s, that contains all
abelian \Cs s and that is closed under direct
sums, inductive limits, and stable isomorphism. Then this class
contains a purely infinite \Cs{} (because it contains all AH-algebras).

\emph{A word of warning:} In the literature, an AH-algebra is often
defined to be an inductive limit of direct sums of building blocks of
the form $pC(\Omega,M_n)p$, where each $\Omega$ is a \emph{compact}
Hausdorff space (and $p$ is a projection in $C(\Omega,M_n)$). With
this definition, AH-algebras always contain non-zero projections. The
algebras we consider, where the building blocs are of the form
$C_0(\Omega,M_n)$ for some \emph{locally compact} Hausdorff space,
should perhaps be called AH$_0$-algebras to distinguish them from the
compact case, but hoping that no confusion will arise, we shall not
distinguish between AH- and AH$_0$-algebras here. 

Every AH-algebra is AF-embeddable. Our Theorem~\ref{thm:A} therefore
gives a new proof of Ozawa's result that there are purely infinite---even
$\cO_\infty$-absorbing---AF-embeddable \Cs s. Moreover, just knowing
that there exists one 
AF-embeddable $\cO_\infty$-absorbing \Cs{}, in combination
with Kirchberg's theorem that all separable, exact \Cs s can be
embedded in $\cO_\infty$, immediately implies that the cone and the
suspension over any separable, exact \Cs{} is AF-embeddable
(Theorem~\ref{thm:B}). This observation yields a new proof of Ozawa's
theorem referred to above.

Section~5 contains some results with relevance 
to the classification program of Elliott. In Section~6 we show that
$\A_\J$ can be embedded into the AF-algebra $\A_\Ca$, where $\Ca$ is
the Cantor set, and hence that the the cone and the
suspension over any separable, exact \Cs{} can be embedded into this
AF-algebra. The ordered
$K_0$-group of $\A_\Ca$ is determined.

I thank Nate Brown for several discussions on quasi-diagonal \Cs s and
on the possibility of embedding quasi-diagonal \Cs s into AF-algebras. I thank
Eberhard Kirchberg for suggesting the nice proof of
Proposition~\ref{prop:AH-AF}, and I thank the referee for suggesting
several improvements to the paper.

\section{The \Cs s $\A_T$} \label{sec:0}

\noindent We review in this section results from Mortensen's paper
\cite{Mor:IO2} on how to associate a \Cs{}
$\A_T$ with each totally ordered, compact, metrizable set $T$, so that
the ideal 
lattice of $\A_T$ is order isomorphic to $T$. Where Mortensen's
algebras are inductive limits of  \Cs s of the form
$C_0(T \setmin \{\max T\},M_{2^n}(\cO_2))$, we consider plain matrix
algebras $M_{2^n}$ in the place  
of $M_{2^n}(\cO_2)$. It turns out that Mortensen's algebras and those
we consider actually are isomorphic when $T = [0,1]$ (see the
second paragraph of Section~5).

Any totally ordered set, which is compact and metrizable in its order
topology, is order isomorphic to a compact subset of $\R$ (where
subsets of $\R$ are given the order structure inherited from $\R$). We shall
therefore assume that we are given a compact subset $T$ of $\R$. 

Put $t_{\max} = \max T$, $t_{\min} = \min T$, and put $T_0 = T \setmin
\{t_{\max}\}$.
Choose a sequence $\{t_n\}_{n=1}^\infty$ in $T_0$ such that the tail
$\{t_k,t_{k+1},t_{k+2}, \dots\}$ is dense in $T_0$ for every $k \in
\N$. Let $\A_T$ be the inductive limit of the sequence
\begin{equation} \label{eq:A}
\xymatrix{C_0(T_0,M_2) \ar[r]^-{\varphi_1} & C_0(T_0,M_4)
  \ar[r]^-{\varphi_2} & C_0(T_0,M_8) \ar[r]^-{\varphi_3} & \cdots
  \ar[r] & \A_T,} 
\end{equation}
where 
\begin{equation} \label{eq:varphi}
\varphi_n(f)(t) = \begin{pmatrix} f(t) & 0 \\ 0 & f(\max\{t,t_n\})
\end{pmatrix} = \begin{pmatrix} f(t) & 0 \\ 0 & (f \circ \chi_{t_n})(t)
\end{pmatrix},
\end{equation}
and where we for each $s$ in $T$ let $\chi_s \colon T \to T$
be the continuous function given by $\chi_s(t) = \max\{t,s\}$. The
algebra $\A_T$ depends a priori on the choice of the dense sequence
$\{t_n\}$. The isomorphism class of $\A_T$ does not depend on this
choice when $T$ is the 
Cantor set (as shown in Section~6) or when $T$ is the interval $\J$
(as will be shown in a forthcoming paper, \cite{KirRor:pi3}). It is
likely that $\A_T$ is independent on $\{t_n\}$ for arbitrary $T$. 

For the sake of brevity, put $A_n = C_0(T_0, M_{2^n}) = C_0(T_0)
\otimes M_{2^n}$. Let $\varphi_{\infty,n} \colon
A_n \to \A_T$ and $\varphi_{m,n} \colon 
A_n \to A_m$, for $n < m$, denote the inductive limit maps, so that
$\A_T$ is the closure of $\bigcup_{n=1}^\infty \varphi_{\infty,n}(A_n)$.

Use the identity $\chi_s \circ \chi_t = \chi_{\max\{s,t\}}$ to see
that 
%\begin{equation} \label{eq:varphi2}
%\varphi_{n+k,n}(f) =  \begin{pmatrix} f_1(t) & 0 & \cdots & 0 \\
%0 & f_2(t) & \cdots & 0 \\ \vdots & \vdots & \ddots & \vdots \\ 0 & 0
%& \cdots & f_{2^k}(t) \end{pmatrix},
%\end{equation}
\begin{equation} \label{eq:varphi2}
\varphi_{n+k,n}(f) =  \begin{pmatrix} f \circ \chi_{\max F_1} & 0 &
  \cdots & 0 \\ 
0 & f \circ \chi_{\max F_2} & \cdots & 0 \\ \vdots & \vdots & \ddots &
\vdots \\ 0 & 0 
& \cdots & f \circ \chi_{\max F_{2^k}} \end{pmatrix},
\end{equation}
%for all natural numbers $n$ and $k$, where
%\begin{equation} \label{eq:varphi3}
%\{f_1,f_2, \dots, f_{2^k}\} = \big\{f \circ \chi_{\max F} \mid F
%\subseteq \{t_n,t_{n+1}, \dots, t_{n+k-1}\} \big\},
%\end{equation}
%(with the convention $\max \emptyset = 0$). 
(with the convention $\max \emptyset = t_{\min} $), where $F_1,F_2,
\dots, F_{2^k}$ is an enumeration of the subsets of $\{t_n,t_{n+1}, 
\dots, t_{n+k-1}\}$. Note
that $\chi_{t_{\min} }$ is the identity map on $T$. 

For each $t \in T$ and for each $n \in \N$ consider the closed ideal 
\begin{equation} \label{eq:I_tn}
I_t^{(n)} \eqdef \{ f \in A_n \mid f(s) = 0 \; \text{when} \; s \ge t \}
\cong C_0(T \cap [t_{\min},t),M_{2^n})
\end{equation}
of $A_n$. Observe that $I_{t_{\min}}^{(n)}=\{0\}$, 
$I_{t_{\max}}^{(n)}=A_n$, and $I_t^{(n)} 
\subset I_s^{(n)}$ whenever $t < s$ for all $n \in \N$. We have
$\varphi_n^{-1}(I_t^{(n+1)}) = I_t^{(n)}$ for all $t$ and for all
$n$, and so
\begin{equation} \label{eq:I_t}
I_t \eqdef \overline{\bigcup_{n=1}^\infty
  \varphi_{\infty,n}(I_t^{(n)})}, \qquad t \in T,
\end{equation}
is a closed two-sided ideal in $\A_T$ such that $I_t^{(n)} =
\varphi_{\infty,n}^{-1}(I_t)$. Moreover, $I_{t_{\min}} = \{0\}$, 
$I_{t_{\max}} = \A_T$, and $I_t \subset I_s$ whenever $s,t \in
T$ and $t < s$.  

\begin{proposition}[cf.\ Mortensen, {\cite[Theorem 1.2.1]{Mor:IO2}}]
\label{prop:ideal-lattice}
Let $T$ be a compact subset of $\R$. Then each closed two-sided ideal in
$\A_T$ is equal to $I_t$ for some $t \in 
T$. It follows that the map $t \mapsto I_t$ is an
order isomorphism from the ordered set $T$ onto the ideal lattice of
$\A_T$. 
\end{proposition}

\begin{proof} Let $I$ be a closed two-sided ideal in
  $\A_T$. Put $I^{(n)} 
  = \varphi_{\infty,n}^{-1}(I)  \vartriangleleft C_0(T_0,M_{2^n}) =
  A_n$, and put 
$$T_n = \bigcap_{f \in I^{(n)}} f^{-1}(\{0\}) \subseteq T, \qquad n
\in \N.
$$
Then $I^{(n)}$ is equal to the set of all continuous functions $f \colon T
\to M_{2^n}$ that vanish on $T_n$. It therefore suffices to show that
there is $t$ in $T$ such that $T_n = T \cap [t,t_{\max}]$ for all
$n$, cf.\ \eqref{eq:I_tn} and \eqref{eq:I_t}. Now,
\begin{equation} \label{eq:T_n}
T_n = T_{n+1} \cup \chi_{t_n}(T_{n+1}) = \bigcup_{F \subseteq X_{n,k}}
\chi_{\max F}(T_{n+k}), \qquad n,k \in \N, 
\end{equation}
where $X_{n,k} = \{t_n,t_{n+1}, \dots, t_{n+k-1}\}$; because if we let
$T'_{n,k}$ denote the right-hand side of \eqref{eq:T_n}, then for all
$f \in C_0(T_0,M_{2^n}) = A_n$,
\begin{eqnarray*}
f|_{T_n} \equiv 0 & \iff &
f \in I^{(n)} \; \iff \; \varphi_{n+k,n}(f) \in I^{(n+k)} \\
& \iff & \forall s \in T_{n+k}: \varphi_{n+k,n}(f)(s) = 0 \\ & 
\overset{\eqref{eq:varphi2}}{\iff} & 
\forall F \subseteq X_{n,k} \; \forall s \in T_{n+k} : f(\chi_{\max
  F}(s)) = 0 \\ & \iff &
f|_{T'_{n,k}} \equiv 0. 
\end{eqnarray*}
It follows from \eqref{eq:T_n} that $\min T_n \le \min T_{n+1}$ for
all $n$; and as  
$$\min \chi_{t_n}(T_{n+1}) = \max\{t_n,\min T_{n+1}\} \ge \min
T_{n+1},$$ 
we
actually have $\min T_n = \min T_{n+1}$ for all $n$. Let $t \in T$ be the
common minimum. Because $t$ belongs to $T_{n+k}$ for all $k$, we can
use \eqref{eq:T_n} to conclude that $T_n$ contains the set
$\{t_n,t_{n+1},t_{n+2}, \dots \} \cap [t,t_{\max}]$; and this set is
by assumption dense in $T \cap (t,t_{\max}]$. This proves the desired
identity: $T_n = T \cap [t,t_{\max}]$, because $T_n$ is a closed subset of
$T \cap [t,t_{\max}]$ and $t$ belongs to $T_n$. 
\end{proof}

\begin{proposition} \label{prop:stable0} 
$\A_T$ is stable for every compact subset $T$ of $\R$.
\end{proposition} 

\begin{proof}
Let $f$ be a positive element in the dense subset $C_c(T_0,M_{2^n})$
of $A_n$ and let $m > n$ be chosen such that $f(t)=0$ for all $t \ge
t_{m-1}$. Then $f \circ \chi_{\max F}
= 0$ for every subset $F$ of $\{t_n,t_{n+1}, \dots, t_{m-1}\}$ that
contains $t_{m-1}$. In the description of $\varphi_{m,n}(f)$
in \eqref{eq:varphi2} we see that $f \circ \chi_{\max F_j} = 0$ for at
least every other $j$. We can therefore find a positive 
function $g$ in $A_{m} = C_0(T_0,M_{2^m})$ such that $g \perp
\varphi_{m,n}(f)$ and $g \sim \varphi_{m,n}(f)$ (the latter in the
sense that $x^*x = g$ and $xx^* = \varphi_{m,n}(f)$ for some $x \in
A_m)$. It follows from \cite[Theorem~2.1 and
Proposition~2.2]{HjeRor:stable} that $\A_T$ is stable 
\end{proof}

\section{A purely infinite AH-algebra}
\label{sec:1} 

\noindent We show in this section that the \Cs{} $\A_\J$ is traceless
and that $\cB = \A_\J \otimes
M_{2^\infty}$ is purely infinite. (In Section~5 it 
will be shown that $\A_\J \cong \cB$.) 

Following \cite[Definition 4.2]{KirRor:pi2} we say that an
exact \Cs{} is \emph{traceless} if it admits no non-zero lower
semi-continuous trace (whose
domain is allowed to be any algebraic ideal of the \Cs). (By
restricting to the case of exact \Cs s we can avoid talking about
quasi-traces; cf.\ Haagerup \cite{Haa:quasi} and Kirchberg
\cite{Kir:quasitraces}.)  

If $\tau$ is a trace defined on  an algebraic ideal $\mathcal{I}$ of a
\Cs{} $B$, and if $I$ is the closure of $\mathcal{I}$, then
$\mathcal{I}$ contains the Pedersen ideal of $I$. In particular,
$(a-\ep)_+$ belongs to $\mathcal{I}$ for every positive element $a$ in
$I$ and for every $\ep>0$. (Here, $(a-\ep)_+ = f_\ep(a)$, where
$f_\ep(t) = \max\{t-\ep,0\}$. Note that $\|a-(a-\ep)_+\| \le
\ep$.) 

\begin{proposition} \label{prop:notraces}
The \Cs{} $\A_\J$ is traceless. 
\end{proposition}

\begin{proof}
Assume, to reach a contradiction, that $\tau$ is a non-zero, lower
semi-continuous, positive trace
defined on an algebraic ideal $\mathcal{I}$ of $\A_\J$, and let $I_t$ be
the closure of $\mathcal{I}$, cf.\
Proposition~\ref{prop:ideal-lattice}. Since $\tau$ is non-zero, $I_t$
is non-zero, and hence $t >0$. 

Identify $I_t^{(n)} = \varphi_{\infty,n}^{-1}(I_t)$ with
$C_0([0,t),M_{2^n})$.  
Put ${\mathcal{I}}^{(n)} = \varphi_{\infty,n}^{-1}({\mathcal{I}})$.
If $x$ is a positive element in $I_t^{(n)}$ and if $\ep >0$, then 
$$\varphi_{\infty,n}((x-\ep)_+) = \big( \varphi_{\infty,n}(x)-\ep
\big)_+ \in \mathcal{I},$$
and so $(x-\ep)_+ \in {\mathcal{I}}^{(n)}$. This shows that
${\mathcal{I}}^{(n)}$ is a dense ideal in $I_t^{(n)}$, and 
hence that ${\mathcal{I}}^{(n)}$ contains $C_c([0,t),M_{2^n})$.

Let $\tau_n$ be the trace on
${\mathcal{I}}^{(n)}$ defined by 
$\tau_n(f) = \tau(\varphi_{\infty,n}(f))$. We show that
\begin{equation} \label{eq:tau_n}
\tau_n(f) = \int_0^t \Tr(f(s)) \, d\mu_n(s), \qquad f \in
C_c([0,t),M_{2^n}),
\end{equation}
for some Radon measure $\mu_n$ on $[0,t)$ (where $\Tr$ denotes the
standard unnormalized trace on $M_{2^n}$). Use 
Riesz' representation theorem to find a Radon measure $\mu_n$ on
$[0,t)$ such that $\tau_n(f) = 2^n\int_0^t f(s) d\mu_n(s)$ for all $f$ in
$C_c([0,t),\C) 
\subseteq C_c([0,t),M_{2^n})$. Let $E \colon C_c([0,t),M_{2^n}) \to
C_c([0,t),\C)$ 
be the conditional expectation given by $E(f)(t) = 2^{-n}\Tr(f(t))$. Then
\begin{equation} \label{eq:E(f)}
E(f) \in \overline{\mathrm{co}}\{ufu^* \mid u \; \text{is a unitary
  element in} \; C([0,t],M_{2^n})\}, \quad f \in C_c([0,t),M_{2^n}), 
\end{equation}
from which we see that $\tau_n(f) = \tau_n(E(f))$. This proves that
\eqref{eq:tau_n} holds. Because 
$\mu_n$ is a Radon measure, $\mu_n([0,s]) < \infty$ for all $s
\in [0,t)$ and for all $n \in \N$. 

Let $\{t_n\}_{n=1}^\infty$ be the sequence in $T$ used in the
definition of $\A_T$.
For each $n$ and $k$ in $\N$ we have $\tau_n = \tau_{n+k} \circ
\varphi_{n+k,n}$. Set $X_{k,n} = \{t_n,t_{n+1}, \dots, t_{n+k-1}\}$
and use \eqref{eq:varphi2} and \eqref{eq:tau_n} to see that
\begin{eqnarray*}
\int_0^t \Tr(f(s)) \, d\mu_n(s) & = & \tau_n(f) \; = \; \tau_{n+k}(
\varphi_{n+k,n}(f)) \\
&=& \int_0^t
\Tr\big(\varphi_{n+k,n}(f)(s)\big) \, d\mu_{n+k}(s) \\
& = & \sum_{F \subseteq X_{k,n}} \int_0^t
\Tr\big((f \circ \chi_{\max(F)})(s)\big) \, d\mu_{n+k}(s) \\
& = & \sum_{F \subseteq X_{k,n}} \int_0^t
\Tr(f(s)) \, d(\mu_{n+k} \circ \chi_{\max(F)}^{-1})(s) 
\end{eqnarray*}
for all $f \in C_c([0,t),M_{2^n})$. This entails that
\begin{equation} \label{eq:mu}
\mu_n =  \sum_{F \subseteq X_{k,n}} \mu_{n+k} \circ
\chi_{\max(F)}^{-1},
\end{equation}
for all natural numbers $n$ and $k$.

We prove next that $\mu_n([0,s])=0$ for all natural numbers $n$ and
for all $s$ in $[0,t)$. Choose $r$ such that $0 < s < r < t$. Put
$Y_{k,n} = X_{k,n} \cap 
[0,s]$ and put $Z_{k,n} = X_{k,n} \cap 
[0,r]$. Observe that
\begin{equation} \label{eq:chi}
\chi_u^{-1}([0,v]) = \begin{cases} \emptyset, & \text{if} \; v <
  u, \\ \mbox{$[0,v]$}, & \text{if} \; v \ge u, \end{cases}
\end{equation}
whenever $u,v \in \J$. Use \eqref{eq:mu} and \eqref{eq:chi} to obtain
\begin{eqnarray} \label{eq:mu0}
\mu_n([0,r]) & = & \sum_{F \subseteq Z_{k,n}} \mu_{n+k}([0,r]) \; = \;
2^{|Z_{k,n}|} \mu_{n+k}([0,r]).
\end{eqnarray}
Use \eqref{eq:mu}, \eqref{eq:chi}, and \eqref{eq:mu0} to see that
\begin{eqnarray*}
\mu_n([0,s]) & = & \sum_{F \subseteq Y_{k,n}} \mu_{n+k}([0,s]) \; = \;
2^{|Y_{k,n}|} \mu_{n+k}([0,s]) \\
& \le & 2^{|Y_{k,n}|} \mu_{n+k}([0,r]) \;
 = \; 2^{-(|Z_{k,n}| - |Y_{k,n}|)} \mu_{n}([0,r]).
\end{eqnarray*}
As
$$\lim_{k \to \infty} \big( |Z_{k,n}| - |Y_{k,n}| \big) = \lim_{k \to \infty}
\big|X_{k,n} \cap \, (s,r] \big| = \infty,$$
(because $\bigcup_{k=n}^\infty X_{k,n} = \{t_n,t_{n+1}, \dots \}$ is
dense in \I), and as 
$\mu_n([0,r]) < \infty$, we conclude that $\mu_n([0,s]) = 0$. It follows
that $\mu_n([0,t)) = 0$, whence $\mu_n$ and $\tau_n$ are
zero for all $n$. 

However, if $\tau$ is non-zero, then $\tau_n$ must be non-zero for
some $n$. To see this, take a positive element $e$ in $\mathcal{I}$
such that $\tau(e)>0$. Because $\tau$ is
lower semi-continuous there is $\ep>0$ such that $\tau((e-\ep)_+) >
0$. Now, ${\mathcal{I}}^{(n)}$ is dense in $I_t^{(n)}$ and
$\bigcup_{n=1}^\infty \varphi_{\infty,n}(I_t^{(n)})$ is dense in
$I_t \supset \mathcal{I}$. It follows that we can 
find $n \in \N$ and a positive element $f$ in ${\mathcal{I}}^{(n)}$ such that
$\|\varphi_{\infty,n}(f) - e\| < \ep$. Use for example \cite[Lemma
2.2]{KirRor:pi2} to find a contraction $d \in A$ such that $d^*
\varphi_{\infty,n}(f)d = (e-\ep)_+$. Put $x =
\varphi_{\infty,n}(f)^{1/2}d$. Then 
\begin{eqnarray*}
\tau_n(f) & = & \tau(\varphi_{\infty,n}(f)) \; \ge \; 
\tau\big(\varphi_{\infty,n}(f)^{1/2}dd^*\varphi_{\infty,n}(f)^{1/2}\big)
\\ & = & \tau(xx^*) \; = \; \tau(x^*x) \; = \; \tau\big((e-\ep)_+\big)
\; > \; 0,
\end{eqnarray*}
and this shows that $\tau_n$ is non-zero.
\end{proof}

\noindent In the formulation of the main result below, $M_{2^\infty}$
denotes the CAR-algebra, or equivalently the UHF-algebra of type
$2^\infty$. 

It is shown in
\cite[Corollary~9.3]{KirRor:pi2} that the following three conditions
are equivalent for a separable, stable (or unital), nuclear \Cs{} $B$:
\begin{enumerate}
\item $B \cong B \otimes \cO_\infty$.
\item $B$ is purely infinite and approximately divisible.
\item $B$ is traceless and approximately divisible.
\end{enumerate}
The \Cs{} $\cO_\infty$ is the Cuntz algebra generated by a sequence
$\{s_n\}_{n=1}^\infty$ of isometries with pairwise orthogonal range
projections. Pure infiniteness of (non-simple) \Cs s was defined in
\cite{KirRor:pi} (see also the introduction). A (possibly non-unital) \Cs{}
$B$ is said to be \emph{approximately divisible} 
if for each natural number $k$ there is a sequence of unital \sh s 
$$\psi_n \colon M_k \oplus M_{k+1} \to \cM(B)$$
such that $\psi_n(x)b-b\psi_n(x) \to 0$ for all $x \in M_k \oplus
M_{k+1}$ and for all $b \in B$, cf.\
\cite[Defini\-tion~5.5]{KirRor:pi}. The tensor product $A \otimes
M_{2^\infty}$ is approximately divisible for any \Cs{} $A$. 

\begin{theorem} \label{thm:A}
Put $\cB = \A_\J \otimes M_{2^\infty}$, where $\A_\J$ is as defined in
\eqref{eq:A}. Then:  
\begin{enumerate}
\item $\cB$ is an inductive limit
$$C_0(\I,M_{k_1}) \to C_0(\I,M_{k_2}) \to C_0(\I,M_{k_3}) \to \cdots
\to \cB,$$
for some natural numbers $k_1,k_2,k_3, \dots$. In particular, $\cB$ is
an AH-algebra. 
\item $\cB$ is traceless, purely infinite, and $\cB \cong \cB \otimes
  \cO_\infty$. 
\end{enumerate}
\end{theorem}

\noindent It is shown in Proposition~\ref{prop:A=B} below that $\A_\J
\cong \cB$. We stress that this fact will not be used in the proof of
Theorem~\ref{thm:B} below.

\begin{proof}
Part (i) follows immediately from the construction of $\A_\J$ and from the fact
that $M_{2^\infty}$ is an inductive limit of matrix algebras.   

(ii). The property of being traceless is preserved after tensoring
with $M_{2^\infty}$, so $\cB$ is traceless by
Proposition~\ref{prop:notraces}. As remarked above, $\cB$ is 
approximately divisible, $\A_\J$ and hence $\cB$ are stable by
Proposition~\ref{prop:stable0}, and as $\cB$ is also nuclear and
separable it follows from \cite[Corollary~9.3]{KirRor:pi2} (quoted above)
that $\cB$ is purely infinite and $\cO_\infty$-absorbing.
\end{proof}

\noindent The \Cs{} $\cB$ is stably projectionless, and, in fact, every
purely infinite  
AH-algebra is (stably) projectionless. Indeed, any projection in
an AH-algebra is finite (in 
the sense of Murray and von Neumann), and any non-zero projection in a
purely infinite \Cs{} is (properly) infinite, cf.\ \cite[Theorem
  4.16]{KirRor:pi}. 

It is impossible to find a \emph{simple} purely infinite AH-algebra,
because all simple purely infinite \Cs s contain properly
infinite projections.

\section{An application to AF-embeddability}

\noindent We show here how Theorem~\ref{thm:A} leads to a new proof of
the recent theorem of Ozawa that the 
cone and the suspension over any exact separable \Cs{} are AF-embeddable,
\cite{Oza:AF}. 

It is well-known that any ASH-algebra, hence any AH-algebra, and hence
the \Cs s $\A_\J$ and $\cB$ from Theorem~\ref{thm:A} are AF-embeddable. For
the convenience of the reader we include a proof of 
this fact---the proof we present is due to 
Kirchberg. (An ASH-algebra is a \Cs{} that arises as the
inductive limit of a sequence of \Cs s each of which is a finite
direct sum of basic building blocks: sub-\Cs s of
$M_n(C_0(\Omega))$---where $n$ and $\Omega$ are allowed to vary.) 

An embedding of $\A_\J$ into an explicit AF-algebra is given in
Section~6. 

\begin{proposition}[Folklore] \label{prop:AH-AF}
Every ASH-algebra admits a faithful embedding into an
AF-algebra.
\end{proposition}

\begin{proof} Note first that if $A$ is a sub-\Cs{} of
  $M_n(C_0(\Omega))$, then its enveloping von Neumann algebra $A^{**}$
  is isomorphic to $\bigoplus_{k=1}^n M_k(\cC_k)$ for some (possibly
  trivial) abelian von Neumann
  algebras $\cC_1, \cC_2, \dots, \cC_n$. If $\cC$ is an abelian von
  Neumann algebra and if
  $D$ is a separable sub-\Cs{} of 
  $M_k(\cC)$, then there is a (separable)
  sub-\Cs{} $D_1$ of 
  $M_k(\cC)$ that contains $D$ and such that $D_1 \cong M_k(C(X))$, where $X$
  is a compact Hausdorff space of dimension zero. In particular, $D_1$
  is an AF-algebra. 

To see 
  this, let $D_0$ be the separable \Cs{} generated by $D$ and the
  matrix units of $M_k \subseteq M_k(\cC)$. Then $D_0 =
  M_k(\cD_0)$ for some 
  separable sub-\Cs{} $\cD_0$ of $\cC$. Any separable sub-\Cs{} of a
  (possibly non-separable) \Cs{} of real rank zero is contained in a
  separable sub-\Cs{} of real rank zero. (This is obtained by successively
  adding projections from the bigger \Cs.) Hence $\cD_0$ is contained in
  a separable real rank zero sub-\Cs{} $\cD_1$ of $\cC$. It follows
  from \cite{BroPed:realrank} that  $\cD_1 \cong C(X)$ for some
  zero-dimensional compact Hausdorff space $X$. Hence $D_1 =
  M_k(\cD_1)$ is as desired.

Assume now that $B$ is an ASH-algebra, so that it is an inductive limit
$$
\xymatrix{
B_1 \ar[r]^-{\psi_1} & B_2 \ar[r]^-{\psi_2} &  B_3 \ar[r]^-{\psi_3} &
\cdots \ar[r] & B,}$$
where each $B_n$ is a finite direct sum of sub-\Cs s of 
$M_m(C_0(\Omega))$. Passing to the bi-dual we get a sequence of finite
von Neumann algebras
$$\xymatrix{
B_1^{**} \ar[r]^-{\psi_1^{**}} & B_2^{**} \ar[r]^-{\psi_2^{**}} &
B_3^{**} \ar[r]^-{\psi_3^{**}} & \cdots }.$$
Use the observation from in the first paragraph (now
applied to direct sums of basic 
building blocks) to
find an AF-algebra $D_1$ such that $B_1 \subseteq D_1
\subseteq B_1^{**}$. Use the observation again to find an
AF-algebra $D_2$ such that $C^*(\psi_1^{**}(D_1),B_2)
\subseteq D_2 \subseteq B_2^{**}$. Continue in this way and find, at
the $n$th stage, an AF-algebra $D_n$ such that
$C^*(\psi_{n-1}^{**}(D_{n-1}),B_n) \subseteq D_n \subseteq
B_n^{**}$. It then follows that the inductive limit $D$ of 
$$\xymatrix{
D_1 \ar[r]^-{\psi_1^{**}} & D_2 \ar[r]^-{\psi_2^{**}} &
D_3 \ar[r]^-{\psi_3^{**}} & \cdots \ar[r] & D}$$
is an AF-algebra that contains $B$.
\end{proof}

\begin{theorem}[Ozawa] \label{thm:B}
The cone $CA = C_0(\I,A)$ over
any separable exact \Cs{} $A$ admits a faithful embedding into an AF-algebra.
\end{theorem} 

\begin{proof} 
By a renowned theorem of Kirchberg any separable exact \Cs{} can be
embedded into the Cuntz algebra $\cO_2$ (see \cite{KirPhi:classI}), and
hence into $\cO_\infty$ (the latter because $\cO_2$ can be
embedded---non-unitally---into $\cO_\infty$). It therefore suffices to
show that $C\cO_\infty = C_0(\I) \otimes \cO_\infty$ is AF-embeddable. 
It is clear from the construction of $\cB$ in Theorem~\ref{thm:A}
that $C_0(\I)$ admits an
embedding into the \Cs{} $\cB$. (Actually, one can embed $C_0(\I)$
into any \Cs{} that absorbs $\cO_\infty$.) As $\cB \cong \cB
\otimes \cO_\infty$, we can embed $C\cO_\infty$ into $\cB$. Now, $\cB$ is
an AH-algebra and therefore AF-embeddable, cf.\
Proposition~\ref{prop:AH-AF}, so $C\cO_\infty$ is AF-embeddable.
\end{proof}

\noindent Ozawa used his theorem in combination with a result of
Spielberg to conclude that the class of AF-embeddable \Cs s is closed
under homotopy invariance, and even more: If $A$ is AF-embeddable
and $B$ is homotopically dominated by $A$, then $B$ is AF-embeddable.

The suspension $SA = C_0((0,1),A)$ is a sub-\Cs{} of $CA$,
and so it follows from Theorem~\ref{thm:B} that also the suspension
over any separable exact \Cs{} is AF-embeddable.

No simple AF-algebra contains a purely infinite
sub-\Cs{}. In fact, any AF-algebra, that has a purely infinite
sub-\Cs{}, must have uncountably many ideals:

\begin{proposition} \label{prop:manyideals}
Suppose that $\varphi \colon A \to B$ is an embedding of a purely infinite
\Cs{} $A$ into an AF-algebra $B$. Let $a$ be a non-zero positive
element in $\image(\varphi)$. For each $t$ in $[0,\|a\|]$ let $I_t$ be the
closed two-sided ideal in $B$ generated by $(a-t)_+$. Then the map $t
\mapsto I_{\|a\|-t}$ defines an injective order embedding of the interval
$[0,\|a\|]$ into the ideal lattice of $B$. 
\end{proposition} 

\begin{proof} Since $A$ is traceless (being purely infinite, cf.\
  \cite{KirRor:pi}), $\image(\varphi) \cap \Dom(\tau)
\subseteq \Ker(\tau)$ for every trace $\tau$ on $B$.

Let $0 \le t < s \le \|a\|$ be given. We show that $I_s$ is strictly
contained in $I_t$.
Find a
  projection $p$ in $\overline{(a-t)_+B(a-t)_+}$ such 
  that $\|(a-t)_+ -p(a-t)_+p\| < s-t$.  There is a trace $\tau$, defined on the
  algebraic ideal in $B$ generated by $p$, with
  $\tau(p)=1$. We claim that
$$I_s \subseteq \Ker(\tau) \subset \Dom(\tau) \subseteq I_t,$$
and this will prove the proposition. To see the first inclusion, there
is $d$ in $B$ such that $(a-s)_+ = d^*p(a-t)_+pd$ (use for example
\cite[Lemma~2.2 and (2.1)]{KirRor:pi2}). Therefore $(a-s)_+$ belongs
to the algebraic ideal generated by $p$, whence $(a-s)_+ \in \image(\varphi)
\cap \Dom(\tau) \subseteq \Ker(\tau)$. This entails that
$I_s$ is contained in the kernel of $\tau$.   

The strict middle inclusion holds because $0 < \tau(p) < \infty$.
The last inclusion holds because
$p$ belongs to $\overline{(a-t)_+B(a-t)_+} \subseteq I_t$. 
\end{proof}

\noindent It follows from Proposition~\ref{prop:stable} below that no
AF-algebra can have ideal lattice isomorphic to 
$[0,1]$, and so the order embedding from
Proposition~\ref{prop:manyideals} can never be surjective.
In Section~6 we show that one can embed a (stably projectionless)
purely infinite \Cs{} into the AF-algebra $\A_\Ca$, where $\Ca$ is the
Cantor set. The ideal lattice of $\A_\Ca$ is the totally ordered and
totally disconnected set $\Ca$. 

\section{Further properties of the algebras $\A_T$}

\noindent Nuclear separable \Cs s that absorb $\cO_\infty$ have
been classified by Kirchberg in terms of an ideal preserving version
of Kasparov's $KK$-theory, see \cite{Kir:fields}. It is
not easy to decide when two such \Cs s with the same primitive ideal 
space are $KK$-equivalent in this sense. There is however a
particularly well understood special case: If $A$ and $B$ are nuclear,
separable, stable \Cs s that absorb the Cuntz algebra $\cO_2$, then
$A$ is isomorphic to $B$ if and only if $A$ and $B$ have homeomorphic
primitive ideal spaces (cf.\ Kirchberg, \cite{Kir:fields}). 

We show in this section that $\A_\J \cong \A_\J \otimes \cO_\infty$ and that
$\A_\J$ is isomorphic to the \Cs{} $\cB$ from Theorem~\ref{thm:A}. It
is shown in a forthcoming paper, \cite{KirRor:pi3}, that $\A_\J \cong
\A_\J \otimes \cO_2$ (using an observartion that $\A_\J$ is zero
homotopic in an ideal-system preserving way, i.e., there is a
$^*$-homomorphism $\Psi \colon A_\J \to C_0([0,1),A_\J)$ such that
$\mathrm{ev}_0 
\circ \Psi = \Id_{\A_\J}$ and $\Psi(J) \subseteq C_0([0,1),J)$ for every
closed two-sided ideal $J$ in $\A_\J$). Thus it follows from
Kirchberg's theorem that $\A_\J$ is the
\emph{unique} separable, nuclear, stable, $\cO_2$-absorbing \Cs{} whose ideal
lattice is (order isomorphic to) $[0,1]$. It seems likely (but remains
open) that any
separable, nuclear, traceless \Cs{} with ideal lattice isomorphic to
$[0,1]$ must absorb $\cO_2$ and hence be isomorphic to $\A_\J$. 

Not all nuclear, separable \Cs s, whose ideal lattice is isomorphic to
$[0,1]$, are purely infinite (or traceless) as shown in
Proposition~\ref{prop:notpi} below.

We derive below a couple of facts about  \Cs s that have
ideal lattice isomorphic to $[0,1]$:

\begin{proposition} \label{prop:stable}
Let $D$ be a \Cs{} with ideal lattice order isomorphic to
$[0,1]$. Then $D$ stably projectionless. If $D$ moreover is
purely infinite and separable, then $D$ is necessarily stable.
\end{proposition}

\begin{proof}
Since $D$ and $D \otimes \cK$ have the same ideal lattice it suffices
to show that $D$ contains no non-zero projections. 
Let $\{I_t \mid t \in [0,1]\}$ be the ideal lattice of $D$ (such that
$I_t \subset I_s$ whenever $t < s$). Suppose, to reach a
contradiction, that $D$ contains a non-zero projection $e$. Let $I_s$
be the ideal in $D$ generated by $e$. The ideal lattice of the unital
\Cs{} $eDe$ is then $\{eI_te \mid t \in [0,s]\}$ and $eI_te \subset
eI_re$ whenever $0 \le t < r \le s$. This is in contradiction with the
well-known fact that any unital \Cs{} has a maximal proper ideal.

Suppose now that $D$ is purely infinite and separable. To show
that $D$ is stable it suffices to show that $D$ has no (non-zero)
unital quotient, cf.\ \cite[Theorem 4.24]{KirRor:pi}. Now, the ideal
lattice of an arbitrary quotient $D/I_s$ of $D$ is equal to $\{I_t/I_s
\mid t \in [s,1]\}$, and this lattice is order isomorphic to the
interval $[0,1]$ (provided that $I_s \ne I_1=D$). It therefore follows
from the first part of the proposition that $D/I_s$ has no non-zero
projection and $D/I_s$ is therefore in particular non-unital. 
\end{proof}

\begin{proposition} \label{prop:A=B}
$\A_T \cong \A_T \otimes M_{2^\infty} \otimes \cK$ for
every compact subset $T$ of $\R$.
\end{proposition}

\begin{proof} It was shown in Proposition~\ref{prop:stable0} that
  $\A_T$ is stable. We proceed to show that $\A_T$ 
  is isomorphic to
  $\A_T \otimes 
  M_{2^\infty}$. Recall that $A_n = C_0(T_0,M_{2^n})$, put
  $\widetilde{A}_n = C(T,M_{2^n})$, and consider the commutative diagram:
$$
\xymatrix{A_1 \ar[r]^-{\varphi_1} \ar@{ )->}[d] & A_2
  \ar[r]^-{\varphi_2}  \ar@{ )->}[d] & A_3 
  \ar[r]^-{\varphi_3} \ar@{ )->}[d] & \cdots \ar[r] & \A_T
  \ar@{ )->}[d] \\ \widetilde{A}_1 \ar[r]_-{\widetilde{\varphi}_1} &
  \widetilde{A}_2 \ar[r]_-{\widetilde{\varphi}_2} & \widetilde{A}_3 
  \ar[r]_-{\widetilde{\varphi}_3} & \cdots \ar[r] & \widetilde{\A},}
$$
where $\varphi_n$ is as defined in \eqref{eq:varphi}, and where
$\widetilde{\varphi}_n \colon \widetilde{A}_n \to
\widetilde{A}_{n+1}$ is defined using the same recipe as in
\eqref{eq:varphi}. The inductive limit \Cs{} $\widetilde{\A}$ is unital,
each $A_n$ 
is an ideal in $\widetilde{A}_n$, and $\A_T$ is (isomorphic to) an ideal in
$\widetilde{\A}$.  

We show that $\widetilde{\A} \cong \widetilde{\A} \otimes
M_{2^\infty}$. This will imply that $\A_T$ is
isomorphic to an ideal of $\widetilde{\A} \otimes M_{2^\infty}$. Each
ideal in  $\widetilde{\A} \otimes M_{2^\infty}$ is of the form $I \otimes
M_{2^\infty}$ for 
some ideal $I$ in $\widetilde{\A}$. As $ M_{2^\infty} \cong
M_{2^\infty} \otimes M_{2^\infty}$ it will follow that
$\A_T \cong \A_T \otimes M_{2^\infty}$. 

By \cite[Proposition~2.12]{BlaKumRor:apprdiv} (and its proof) to prove
that $\widetilde{\A} \cong \widetilde{\A} \otimes M_{2^\infty}$ 
it suffices to show that 
for each finite subset $G$ of $\widetilde{\A}$ and for each $\ep>0$ there is a
unital \sh{} $\lambda \colon M_2 \to \widetilde{\A}$ such that
$\|\lambda(x)g-g\lambda(x)\| \le \ep\|x\|$ for all $x \in M_2$ and for
all $g \in G$. We may assume that $G$ is contained in
$\widetilde{\varphi}_{\infty,n}(\widetilde{A}_n)$ for some natural
number $n$. Put 
$H= \widetilde{\varphi}_{\infty,n}^{-1}(G)$. It now suffices to find a
natural number $k$ and a unital \sh{} $\lambda \colon M_2 \to
\widetilde{A}_{n+k}$ such that 
\begin{equation} \label{eq:b}
\|\lambda(x)\widetilde{\varphi}_{n+k,n}(h) -
\widetilde{\varphi}_{n+k,n}(h)\lambda(x)\| \le \ep\|x\|, \qquad x \in
M_2, \quad h \in H.
\end{equation}

Put $t_{\min} = \min T$, and find $\delta>0$ such that
$\|h(t)-h(t_{\min})\| \le \ep/2$ for all $h$ in $H$ and for all $t$ in $T$
with $|t-t_{\min}| < \delta$. Let $\{t_n\}$ be the dense sequence in
$T_0$ used in the definition of $\A_T$.
Find $m \ge n$ such that $|t_m-t_{\min}| < \delta$. Put $k =
m+1-n$, and organize the elements in $X= \{t_n,t_{n+1},
\dots, t_{m+1}\}$ in increasing order and relabel the elements by $s_1 \le
s_2 \le s_3 \le \dots \le s_k$. Let $F_1,F_2,
\dots, F_{2^k}$ be the subsets of $X$ ordered such that
$F_1=\emptyset$ and
$$\max F_2 = s_1, \; \max F_3 = \max F_4 = s_2,
\; \dots, \; \max F_{2^{k-1}+1} = \cdots = \max F_{2^k} = s_k.$$
Then $|s_1-t_{\min}| < \delta$, and so $\|h \circ \chi_{\max F_{1}} - h \circ
\chi_{\max F_{2}}\| \le \ep$ for $h \in H$ (we use the convention
$\max \emptyset = t_{\min}$); and  
$h \circ \chi_{\max F_{2j-1}} = h \circ \chi_{\max F_{2j}}$ when
$j \ge 2$ for all $h$.  

We shall use the picture of
$\varphi_{n+k,n}$ given in \eqref{eq:varphi2}, which is valid also for
$\widetilde{\varphi}_{n+k,n}$. However, since the sets $F_1, F_2,
\dots, F_k$ possibly have been 
permuted, $\varphi_{n+k,n}$ and the expression in \eqref{eq:varphi2}
agree only up to unitary equivalence.
Let $\lambda \colon M_2 \to \widetilde{A}_{n+k}$ be the unital \sh{}
given by $\lambda(x) = \diag(x,x, \dots, x)$ (with 
$2^{k-1}$ copies of $x$). Use \eqref{eq:varphi2} and
the estimate 
\begin{eqnarray*}
 & & \big\|x \begin{pmatrix} h \circ \chi_{\max F_{2j-1}} & 0 \\ 0 & h \circ
  \chi_{\max F_{2j}}  \end{pmatrix} - 
 \begin{pmatrix} h \circ \chi_{\max F_{2j-1}} & 0 \\ 0 & h \circ
  \chi_{\max F_{2j}}  \end{pmatrix} x\big\|  \\ & \le &
\|x\|\|h \circ \chi_{\max {F_{2j-1}}}-h \circ \chi_{\max {F_{2j}}}\| \
\le \; \ep\|x\|,
\end{eqnarray*}
that holds for $j=1,2, \dots, 2^{k-1}$, for $h \in H$, and for all $x \in
M_2(\C) \subseteq C(T,M_2)$, to conclude that \eqref{eq:b} holds, and
hence that $\widetilde{A} \cong \widetilde{A} \otimes M_{2^\infty}$.
\end{proof}

\noindent
Propositions~\ref{prop:A=B} together with Theorem~\ref{thm:A}
yield: 

\begin{corollary} \label{cor:A_I}
The \Cs{} $\A_\J$ is purely infinite and $\A_\J \cong
\A_\J \otimes \cO_\infty$. 
\end{corollary}

\noindent We conclude this section by showing that the tracelessness of
the \Cs s $\A_\J$ (established in Proposition~\ref{prop:notraces}) is not
a consequence of its ideal lattice being isomorphic to
$[0,1]$.

%\newpage
\begin{proposition} \label{prop:notpi}
Let $\{l_n\}_{n=1}^\infty$ be a sequence of positive integers, and let
$\{t_n\}_{n=1}^\infty$ be a dense sequence in $\I$. Put $k_1=1$,  
put $k_{n+1} = (l_n+1)k_n$ for $n \ge 1$, and put $D_n =
C_0(\I,M_{k_n})$. Let $\cD$ be the inductive limit of the sequence
$$\xymatrix{D_1 \ar[r]^-{\psi_1} & D_1 \ar[r]^-{\psi_2} & D_2
  \ar[r]^-{\psi_3} & \cdots \ar[r] & \cD,}$$
where $\psi_n(f) = \diag(f,f, \dots, f, f\circ \chi_{t_n})$ (with
$l_n$ copies of $f$), and where $\chi_{t_n} \colon \J \to \J$ as
before is given by $\chi_{t_n}(s) = \max\{s,t_n\}$. 

It follows that the ideal lattice of $\cD$ is isomorphic to
the interval $\J$. Moreover, if $\prod_{n=1}^\infty l_n/(l_n+1) >0$,
then $\cD$ has a non-zero bounded trace, in which case $\cD$ is not
stable and not purely infinite. 
\end{proposition}

\begin{proof} An obvious modification of the proof of
  Proposition~\ref{prop:ideal-lattice} shows that the ideal lattice of
  $\cD$ is isomorphic to $[0,1]$. As in the proof of
  Proposition~\ref{prop:A=B} we construct a unital \Cs{}
  $\widetilde{\cD}$, in which $\cD$ is a closed two-sided ideal, by letting
  $\widetilde{\cD}$ be the inductive limit of the sequence
$$
\xymatrix{D_1 \ar[r]^-{\psi_1} \ar@{ )->}[d] & D_2
  \ar[r]^-{\psi_2}  \ar@{ )->}[d] & D_3 
  \ar[r]^-{\psi_3} \ar@{ )->}[d] & \cdots \ar[r] & \cD
  \ar@{ )->}[d] \\ \widetilde{D}_1 \ar[r]_-{\widetilde{\psi}_1} &
  \widetilde{D}_2 \ar[r]_-{\widetilde{\psi}_2} & \widetilde{D}_3 
  \ar[r]_-{\widetilde{\psi}_3} & \cdots \ar[r] & \widetilde{\cD},}
$$
where $\widetilde{D}_n = C(\J,M_{k_n})$ and $\widetilde{\psi}_n(f) =
\diag(f, \dots, f, f \circ \chi_{t_n})$. Remark that 
$$\widetilde{\cD}/\cD \; \cong \; \lim_{\longrightarrow} \widetilde{D}_n/D_n 
\; \cong \; \lim_{\longrightarrow} M_{k_n}$$
is a UHF-algebra. If $\tau$ is a tracial state on $\widetilde{\cD}$ that
vanishes on $\cD$, then $\tau$ is the composition of the quotient
mapping $\widetilde{\cD} \to \widetilde{\cD}/\cD$ and the unique tracial state
on the UHF-algebra $\widetilde{\cD}/\cD$. It follows that there is only one
tracial state $\tau$ on $\widetilde{\cD}$ that vanishes on $\cD$. 

Suppose now that $\prod_{n=1}^\infty l_n/(l_n+1) >0$. It then follows,
as in the construction of Goodearl in \cite{Goo:goodearl_alg}, that
the simplex of 
tracial states on $\widetilde{\cD}$ is homeomorphic to the simplex of
probability measures on $[0,1]$ and hence that $\widetilde{\cD}$ has a
tracial state that does not vanish on $\cD$. The restriction of this
trace to $\cD$ is then the desired non-zero bounded trace. (Goodearl
constructs simple \Cs s; and where $f \circ \chi_{t_n}$ appears in our
connecting map $\widetilde{\psi}_n$, Goodearl uses a point evaluation,
i.e., the constant function $t \mapsto f(t_n)$. Goodearl's proof
can nonetheless and without changes be applied in our situation.)
\end{proof}

\section{An embedding into a concrete AF-algebra}

\noindent Let $T$ be a compact subset of $\R$ and set $T_0 = T \setmin
\{\max T\}$. Then $C_0(T_0,M_{2^n})$ is an AF-algebra if and only if
$T$ is totally disconnected. It follows that the \Cs{} $\A_T$ (defined
in \eqref{eq:A}) is an AF-algebra whenever $T$ is totally disconnected. Let
$\Ca$ denote the Cantor set (realized as the ``middle third'' subset of
\J, and with the total order it inherits from its embedding in
$\R$). Actually any totally disconnected, compact subset of $\R$ with
no isolated points is order isomorphic to $\Ca$. 

We show here that the AF-algebra from Theorem~\ref{thm:B}, into which
the cone over any separable exact \Cs{} can be embedded, 
can be chosen to be $\A_\Ca$. The ideal lattice of
$\A_\Ca$ is order isomorphic to $\Ca$ (by
Proposition~\ref{prop:ideal-lattice}). In the light of
Proposition~\ref{prop:manyideals} and by the fact that the ideal lattice
of an AF-algebra is totally disconnected (in an appropriate sense)
the AF-algebra $\A_\Ca$ has the least complicated ideal lattice among
AF-algebras that admit embeddings of (stably projectionless) purely
infinite \Cs s.  

We begin by proving a general result on when $\A_{S}$ can be
embedded into $\A_T$:

\begin{proposition} \label{prop:embedding1}
Let $S$ and $T$ be compact subsets of $\R$. Set $T_0 =  T
\setminus \{\max T\}$ and $S_0 = S \setminus \{ \max S \}$. 
Suppose there is a
continuous, increasing, surjective function $\lambda \colon T \to
S$ such that $\lambda(T_0) = S_0$. Let $\{t_n\}_{n=1}^\infty$ be a
sequence in $T_0$ such that 
$\{t_n\}_{n=k}^\infty$ is dense in $T_0$ for every $k$, and put $s_n =
\lambda(t_n)$. Then
$\{s_n\}_{n=k}^\infty$ is dense in $S_0$ for every $k$, and 
there is an 
injective \sh{} $\lambda^\sharp \colon \A_{S} \to \A_T$, when $\A_T$ and
$\A_S$ are inductive limits as in \eqref{eq:A} with respect to the
sequences $\{t_n\}_{n=1}^\infty$ and $\{s_n\}_{n=1}^\infty$, respectively.
If $\lambda$ moreover is injective, then $\lambda^\sharp$ is an isomorphism. 
\end{proposition}

\begin{proof} There is a commutative diagram:
\begin{equation} \label{eq:iota}
%\begin{split}
\xymatrix{C_0(S_0,M_2) \ar[r]^-{\varphi_1} \ar[d]_{\widehat{\lambda}}
  & C_0(S_0,M_4) \ar[r]^-{\varphi_2} \ar[d]_{\widehat{\lambda}}
 & C_0(S_0,M_8) \ar[r]^-{\varphi_3} \ar[d]_{\widehat{\lambda}} & \cdots
 \ar[r] & \A_{S} \ar@{-->}[d]^-{\lambda^\sharp}\\ 
C_0(T_0,M_2) \ar[r]_-{\psi_1}
  & C_0(T_0,M_4) \ar[r]_-{\psi_2}
 & C_0(T_0,M_8) \ar[r]_-{\psi_3} & \cdots
 \ar[r] & \A_T}
%\end{split}
\end{equation}
where $\widehat{\lambda}(f) = f \circ \lambda$, and where
\begin{equation} \label{eq:varphi-psi}
\varphi_n(f) = \begin{pmatrix} f & 0 \\ 0 & f \circ \chi_{s_n}
\end{pmatrix}, \qquad \psi_n(f) = \begin{pmatrix} f & 0 \\ 0 & f \circ
  \chi_{t_n} \end{pmatrix},
\end{equation}
cf.\ \eqref{eq:varphi}. Note that $\lambda(t_{\max}) = s_{\max}$
(because $\lambda$ is surjective), and so  $\widehat{\lambda}(f)(t_{\max}) =
f(\lambda(t_{\max})) = f(s_{\max}) = 0$. 
To see that the diagram \eqref{eq:iota} indeed is commutative we
must check that $\widehat{\lambda} \circ \varphi_n = \psi_n \circ
\widehat{\lambda}$ for all $n$. By \eqref{eq:varphi-psi},
$$(\widehat{\lambda} \circ \varphi_n)(f) =  \begin{pmatrix} f \circ
  \lambda & 0 \\ 0
  & f \circ \chi_{s_n} \circ \lambda \end{pmatrix}, \qquad 
(\psi_n \circ \widehat{\lambda})(f) =  \begin{pmatrix} f \circ
  \lambda & 0 \\ 0
  & f \circ \lambda \circ \chi_{t_n} \end{pmatrix}, \qquad 
$$
for all $f \in C_0(S_0,M_{2^n})$, so it suffices
to check that $\chi_{s_n} \circ \lambda = \lambda \circ
\chi_{t_n}$. But
$$(\chi_{s_n} \circ \lambda)(x) = \max\{\lambda(x),s_n\} =
\max\{\lambda(x),\lambda(t_n)\} = \lambda\big( \max\{x,t_n\} \big) =
(\lambda \circ \chi_{t_n})(x),$$
where the third equality holds because $\lambda$ is
increasing. 

Each map $\widehat{\lambda}$ in the diagram \eqref{eq:iota} is
injective (because $\lambda$ is surjective), so the \sh{} ${\lambda^\sharp}
\colon \A_{S} \to \A_T$ induced by the diagram is injective. 

If $\lambda$ also is injective, then each map $\widehat{\lambda}$ in
\eqref{eq:iota} is an isomorphism in which case ${\lambda^\sharp}$ is an
isomorphism. 
\end{proof}

\noindent Combine (the proof of) Theorem~\ref{thm:B} with
Proposition~\ref{prop:A=B} to obtain:

\begin{proposition} \label{prop:embedding2}
The cone and the suspension over any separable exact \Cs{} admits an
embedding into the AH-algebra $\A_\J$. 
\end{proposition}

\begin{lemma} \label{lm:Cantor} There is a
  continuous, increasing, surjective map $\lambda \colon \Ca \to
  [0,1]$ that maps $[0,1)$ into $\Ca_0$, where $\Ca$ is the Cantor set
  and where $\Ca_0 = \Ca \setmin \{1\}$. 
\end{lemma}

\begin{proof} Each $x$ in $\Ca$ can be written $x =
  \sum_{n \in F} 2 \gange 3^{-n}$ for a unique subset $F$ of $\N$. We
  can therefore define $\lambda$ by 
$$\lambda\big(\sum_{n \in F} 2 \gange 3^{-n} \big) = \sum_{n \in F} 
2^{-n}, \qquad F \subseteq \N.$$
It is straightforward to check that $\lambda$ has the desired properties.
\end{proof} 

\begin{corollary} \label{cor}
The cone and the suspension over any separable exact \Cs{} admits an
embedding into the AF-algebra $\A_\Ca$.  
\end{corollary}

\begin{proof}
It follows from Proposition~\ref{prop:embedding1} and
Lemma~\ref{lm:Cantor} that $\A_\J$ can be embedded into $\A_\Ca$. The
corollary is now an immediate consequence of
Proposition~\ref{prop:embedding2}.  
\end{proof}

\noindent By a renowned theorem of
Elliott, \cite{Ell:AF}, the ordered $K_0$-group is a complete
invariant for the stable isomorphism class of an
AF-algebra. We shall therefore go to some length to calculate the ordered
group $K_0(\A_\Ca)$. 

As $K_0(\A_\Ca)$ does not depend on the choice of dense sequence
$\{t_n\}_{n=1}^\infty$ used in 
the inductive limit description of $\A_\Ca$, \eqref{eq:A}, it follows
in particular from Proposition~\ref{prop:dimgp} below that the isomorphism
class of $\A_\Ca$ is independent of this sequence.

The Cantor set $\Ca$ is realized as the ``middle-third''
subset of $[0,1]$ (so that $0 = 
\min \Ca$ and $1 = \max \Ca$). 
Consider the countable abelian group $G=C_0(\Ca_0,\Z[\frac{1}{2}])$ where the
composition is addition, and where the group of Dyadic rationals
$\Z[\frac{1}{2}]$ is given the discrete topology. Equip $G$ with the
lexicographic 
order, whereby $f \in G^+$ if and only if either $f=0$ or $f(t_0)>0$
for $t_0 = \sup\{t \in \Ca \mid f(t) \ne 0\}$. (The set $\{t \in \Ca
\mid f(t) \ne 0\}$ is clopen because $\Z[\frac{1}{2}]$ is discrete, and so
$f(t_0) \ne 0$.) It is easily checked that $(G,G^+)$ is a totally
ordered abelian group, and hence a dimension group.

\begin{proposition} \label{prop:dimgp}
The group $K_0(\A_\Ca)$ is order isomorphic to the group
$C_0(\Ca_0,\Z[\frac{1}{2}])$ equip\-ped with the lexicographic ordering. 
\end{proposition}

\begin{proof}
Let $\{t_n\}_{n=1}^\infty$ be any
sequence in $\Ca_0 = \Ca \setmin \{1\}$ such that
$\{t_k,t_{k+1},t_{k+2}, \dots\}$ is dense in $\Ca_0$ for all
$k$. Write $\A_\Ca$ as an inductive limit with connecting maps
$\varphi_n$ as in \eqref{eq:A}.

By continuity of $K_0$ and because $K_0(C_0(\Ca_0,M_{2^n})) \cong
C_0(\Ca_0, \Z)$ (as ordered abelian groups) (see eg.\
\cite[Exercise~3.4]{RorLarLau:k-theory}), the ordered abelian group
$K_0(\A_\Ca)$ is the inductive limit of the sequence 
$$\xymatrix{C_0(\Ca_0,\Z) \ar[r]^-{\alpha_1} & C_0(\Ca_0,\Z)
  \ar[r]^-{\alpha_2} & C_0(\Ca_0,\Z) \ar[r]^-{\alpha_3} & \cdots \ar[r]
  & K_0(\A_\Ca),
}
$$
where $\alpha_n(f) = K_0(\varphi_n)(f) = f + f \circ \chi_{t_n}$. 

Choose for each $n \in \N$ a partition $\{A^{(n)}_1, A^{(n)}_2,
\dots, A^{(n)}_{2^n}\}$ of $\Ca$ into clopen intervals (written in
increasing order) such that
\begin{itemize}
\item[(a)] $A^{(n)}_j = A^{(n+1)}_{2j-1} \cup A^{(n+1)}_{2j}$,
\item[(b)] $t_n \in A^{(n)}_1$ for infinitely many $n$,
\item[(c)] $\bigcup_{n=1}^\infty \{A^{(n)}_1, A^{(n)}_2,
\dots, A^{(n)}_{2^n}\}$ is a basis for the topology on $\Ca$.
\end{itemize}
Set
$\cF=\bigcup_{n=1}^\infty \{A^{(n)}_1, A^{(n)}_2, 
\dots, A^{(n)}_{2^n-1}\}$, and set
$$H_n = \Span \{1_{A^{(n)}_j} \mid j = 1,2, \dots, 2^n-1 \} \subseteq
C_0(\Ca_0,\Z).$$
Note that $1_{A^{(n)}_{2^n}}$ does not belong to $C_0(\Ca_0,\Z)$
because $1 \in A^{(n)}_{2^n}$.

We outline the idea of the rather lengthy proof below. We show first that
$\alpha_n(H_n) \subseteq H_{n+1}$ for all $n$ and that
$\bigcup_{n=1}^\infty \alpha_{\infty,n}(H_n) = K_0(\A_\Ca)$, where
$\alpha_{\infty,n} = K_0(\varphi_{\infty,n})$ is the inductive limit
  homomorphism from $C_0(\Ca_0,\Z)$ to $K_0(\A_\Ca)$. We then
  construct positive, injective group homomorphisms $\beta_n \colon H_n \to
  C_0(\Ca_0, \Z[\frac{1}{2}])$ that satisfy $\beta_{n+1} \circ
  \alpha_n = \beta_n$ for all $n$, and which therefore induce a positive
  injective group homomorphism $\beta \colon K_0(\A_\Ca) \to C_0(\Ca_0,
  \Z[\frac{1}{2}])$. It is finally proved that $\beta$ is onto and
  that $K_0(\A_\Ca)$ is totally ordered, and from this one can
  conclude that $\beta$ is an order isomorphism. 

For each interval $[r,s] \cap \Ca$ and for each $t \in \Ca$,
\begin{equation} \label{eq:1_A}
1_{[r,s] \cap \Ca} \circ \chi_t = \begin{cases} 1_{[r,s] \cap \Ca}, &
  t < r,\\ 
1_{[0,s] \cap \Ca}, & r \le t \le s,\\
0, & t > s. \end{cases}
\end{equation}
Suppose that
$A_1,A_2, \dots, A_m$ is a partition of $\Ca$ into clopen intervals,
written in increasing order, and that $t \in A_{j_0}$. Then, by
\eqref{eq:1_A},  
\begin{equation} \label{eq:alpha}
1_{A_j} + 1_{A_j} \circ \chi_t = \begin{cases} 
1_{A_j}, & j < j_0, \\
2 \gange 1_{A_{j}}+1_{A_{j-1}} + \cdots + 1_{A_{1}}, & j = j_0,\\
2 \gange 1_{A_j}, & j >j_0.
\end{cases}
\end{equation}
The lexicographic order on $G =
C_0(\Ca_0,\Z[\frac{1}{2}])$ has the following description: If $k \le n$ and
if $r_1, r_2, \dots, 
r_k$ are elements in $\Z[\frac{1}{2}]$ with $r_k \ne 0$, then
\begin{equation} \label{eq:order}
r_k 1_{A_k} + r_{k-1} 1_{A_{k-1}} + \cdots + r_1 1_{A_1} \in G^+ \iff
r_k > 0. 
\end{equation}

It follows from \eqref{eq:alpha} that $\alpha_n(H_n) = H_n \subseteq
H_{n+1}$. As $\cF$ is a basis for the topology of $\Ca$, the set
$\{1_A \mid A \in \cF\}$ generates $C_0(\Ca_0,\Z)$. 
To prove that $\bigcup_{n=1}^\infty \alpha_{\infty,n}(H_n) =
K_0(\A_\Ca)$ it suffices to show that $\alpha_{\infty,m}(1_A)$ belongs to
$\bigcup_{n=1}^\infty \alpha_{\infty,n}(H_n)$ for every $A$ in $\cF$
and for every $m$ in $\N$. Take $A \in \cF$ and 
find a natural number $n \ge m$ such that $1_A$ belongs to $H_n$. Let
$A'$ be the clopen interval in $\Ca$ consisting of all points in $\Ca$
that are smaller than $\min A$. Then $1_{A'}$ belongs to $H_n$, and
$\alpha_{n,m}(1_A)$ belongs to 
$\Span\{1_{A'},1_A\} \subseteq H_n$ by \eqref{eq:alpha}. Hence
$\alpha_{\infty,m}(1_A) = \alpha_{\infty,n}(\alpha_{n,m}(1_A))$ 
belongs to $\alpha_{\infty,n}(H_n)$.

The next step is to find a sequence of positive, injective group
homomorphisms $\beta_n \colon H_n \to G$ such that $\beta_{n+1} \circ
\alpha_n = \beta_n$. (This sequence will then induce a positive,
injective group homomorphism $\beta \colon K_0(\A_\Ca) \to G$.) Each function 
$\{1_{A_1^{(n)}}, 1_{A_2^{(n)}}, \dots, 1_{A_{2^n-1}^{(n)}}\}
\to G^+$ extends uniquely to a positive group homomorphism $H_n
\to G$, and so it suffices to specify $\beta_n$ on this generating set. We
do so by setting
\begin{equation} \label{eq:beta}
\beta_n(1_{A_j^{(n)}}) = \delta(j,j,n) 1_{A_j^{(n)}} +
\sum_{i=1}^{j-1} \delta(j,i,n) 1_{A_i^{(n)}}, \qquad j = 1,2, \dots, 2^n-1,
\end{equation}
for suitable coefficients, $\delta(j,i,n)$,  in $\Z[\frac{1}{2}]$---to
be contructed---such that
$\delta(j,j,n) = 2^{-k} >0$ for some $k \in \N$, and such that
$1_{A_j^{(n)}}$ belongs to the image of $\beta_n$ for $j=1,2,
\dots, 2^n-1$. Positivity of $\beta_n$ will follow from
\eqref{eq:order}, \eqref{eq:beta}, and the fact that $\delta(j,j,n)>0$.

For $n=1$ set $\beta_1(1_{A_1^{(1)}})=1_{A_1^{(1)}}$, so that
$\delta(1,1,1)=1$. Suppose that $\beta_n$ has been found. The point
$t_n$ belongs to $A_{j_0}^{(n)}$ for some $j_0$. The equation
$\beta_{n+1}(\alpha_n(1_{A_j^{(n)}})) = \beta_n(1_{A_j^{(n)}})$ has by
\eqref{eq:alpha} the solution:
\begin{equation} \label{eq:beta2}
\beta_{n+1}(1_{A_j^{(n)}}) = \begin{cases}
\beta_n(1_{A_j^{(n)}}), & j<j_0, \\
\tfrac{1}{2}\beta_n(1_{A_j^{(n)}}) - \tfrac{1}{2} \sum_{i=1}^{j-1}
\beta_n(1_{A_i^{(n)}}), & j = j_0,\\
\tfrac{1}{2}\beta_n(1_{A_j^{(n)}}), & j > j_0.
\end{cases}
\end{equation}
Extend $\beta_{n+1}$ from $H_n$ to $H_{n+1}$ as follows: 
\begin{eqnarray*}
\beta_{n+1}(1_{A_{2j-1}^{(n+1)}}) & = & \delta(j,j,n) 1_{A_{2j-1}^{(n+1)}} +
\sum_{i=1}^{j-1} \delta(j,i,n) 1_{A_{i}^{(n)}}, \quad j = 1,
\dots, j_0-1,\\ 
\beta_{n+1}(1_{A_{2j}^{(n+1)}}) & = & \delta(j,j,n)
1_{A_{2j}^{(n+1)}},  \quad j = 1,  \dots, j_0-1,\\
\beta_{n+1}(1_{A_{2j-1}^{(n+1)}}) & = & \frac{1}{2}\delta(j,j,n)
1_{A_{2j-1}^{(n+1)}} + \frac{1}{2} \sum_{i=1}^{j-1} \big(
\delta(j,i,n) - \sum_{k=i}^{j-1} \delta(k,i,n)\big)
1_{A_{i}^{(n)}}, \quad j = j_0,\\
\beta_{n+1}(1_{A_{2j}^{(n+1)}}) & = & \frac{1}{2}\delta(j,j,n)
1_{A_{2j}^{(n+1)}}, \quad j = j_0,\\
\beta_{n+1}(1_{A_{2j-1}^{(n+1)}}) & = & \frac{1}{2}\delta(j,j,n)
1_{A_{2j-1}^{(n+1)}} + 
\frac{1}{2}\sum_{i=1}^{j-1} \delta(j,i,n) 1_{A_{i}^{(n)}}, \quad
j = j_0+1, \dots, 2^n-1,\\
\beta_{n+1}(1_{A_{2j}^{(n+1)}}) & = & \frac{1}{2}\delta(j,j,n)
1_{A_{2j}^{(n+1)}},  \quad
j = j_0+1, \dots, 2^n-1,\\
\beta_{n+1}(1_{A_{2^n-1}^{(n+1)}}) & = & 1_{A_{2^n-1}^{(n+1)}}.
\end{eqnarray*}
The coefficients, implicit in these expressions for
$\beta_{n+1}(1_{A_{j}^{(n+1)}})$, will be our
$\delta(j,i,n+1)$. 

It follows by induction on $n$ that each coefficient
$\delta(j,i,n)$ belongs to $\Z[\frac{1}{2}]$ and that $\delta(j,j,n) =
2^{-k}$ for some $k \in \N$ (that depends on $j$ and $n$). The
formula above for $\beta_{n+1}$ is consistent with \eqref{eq:beta2},
and so $\beta_{n+1} \circ \alpha_n = \beta_n$. It also follows by
induction on $n$ that $1_{A_j^{(n)}}$ belongs to
$\image(\beta_{n})$ for $j = 1, 2, \dots, 2^n-1$. This clearly holds
for $n=1$. Assume it holds for some $n \ge 1$. Then $1_{A_j^{(n)}}$
belongs to 
$\image(\beta_{n}) \subseteq \image(\beta_{n+1})$ for $j = 1, 2,
\dots, 2^n-1$, and hence
$1_{A_{2j}^{(n+1)}}$, $1_{A_{2j-1}^{(n+1)}} =
1_{A_j^{(n)}}-1_{A_{2j}^{(n+1)}}$, and $1_{A_{2^n-1}^{(n+1)}}$ belong to
$\image(\beta_{n+1})$. It is now verified that each $\beta_n$ is as desired. 

To complete the proof we must show that the positive, injective, group
homomorphism $\beta \colon K_0(\A_\Ca) \to G$ is surjective and that
$\beta(K_0(\A_\Ca)^+) = G^+$. The former follows from the already
established fact that $1_A$ belongs to the image of $\beta$ for all $A
\in \cF$, and from the fact, which follows from Proposition~\ref{prop:A=B},
that if $f$ belongs to $\image(\beta)$, then so does
$\tfrac{1}{2}f$. The latter identity is proved by verifying that
$K_0(\A_\Ca)$ is totally ordered. 

%Assume that $f$ belongs to $\image(\beta)$. Then $f$ belongs to
%$\image(\beta_n)$ for some $n$. By (b) we can find $m \ge n$ such that
%$t_m \in A_{1}^{(m)}$. Now, $f \in \image(\beta_n) \subseteq
%\image(\beta_m)$ and so $f = \beta_m(g)$ for some $g$ in $H_m$. Use
%\eqref{eq:alpha} to see that $\alpha_m(g) = 2g$, and hence that 
%$$\tfrac{1}{2}f = \tfrac{1}{2} \beta_{m}(g) =
%\tfrac{1}{2} \beta_{m+1}(\alpha_m(g)) =
%\beta_{m+1}(g) \in \image(\beta).$$

To show that $K_0(\A_\Ca)$ is totally ordered we
must show that either $f$ or $-f$ is positive for each non-zero $f$ in
$K_0(\A_\Ca)$. Write $f = \alpha_{\infty,n}(g)$ for a suitable $n$ and
$g \in C_0(\Ca_0,\Z)$. Let $r$ be the largest point in $\Ca$ for
which $g(r) \ne 0$. Upon replacing $f$ by $-f$, if necessary, we can
assume that $g(r)$ is positive. There is a (non-empty) clopen
interval $A=[s,r] \cap \Ca$ for which $g(t) \ge 1$ for all $t$ in
$A$. Put $X_{k,n} = \{t_n,t_{n+1}, \dots, t_{n+k-1}\}$, $Y_{k,n} =
X_{k,n} \cap [0,r]$, and $Z_{k,n} = X_{k,n} \cap [0,s)$. By
\eqref{eq:1_A} and an analog of  \eqref{eq:varphi2} we get 
\begin{eqnarray*} 
\alpha_{n+k,n}(g)  & =  & \sum_{F \subseteq X_{k,n}} g \circ
\chi_{\max F}  \; = \; \sum_{F \subseteq Y_{k,n}} g \circ \chi_{\max
  F} \\ & \ge  & \sum_{F \subseteq Z_{k,n}} \min g(\Ca_0) +
\sum_{F \subseteq Y_{k,n}, \, F \nsubseteq Z_{k,n}} 1_A \circ  \chi_{\max F} 
\\ & = &  2^{|Z_{k,n}|} \gange \min g(\Ca_0) + \big(2^{|Y_{k,n}|} -
 2^{|Z_{k,n}|} \big) \gange 1_{[0,r] \cap \Ca}.  
\end{eqnarray*}
Now,
$$\lim_{k \to \infty} \big(|Y_{k,n}| - |Z_{k,n}|\big) = \lim_{k \to \infty}
|X_{k,n} \cap [r,s]| = \infty,$$ 
so $\alpha_{n+k,n}(g) \ge 0$ for some large enough $k$. But then $f =
\alpha_{\infty,n+k}(\alpha_{n+k,n}(g))$ is positive.
\end{proof}

%\newpage

{\small{
\bibliographystyle{amsplain}
%\bibliography{operator}
\providecommand{\bysame}{\leavevmode\hbox to3em{\hrulefill}\thinspace}
\providecommand{\MR}{\relax\ifhmode\unskip\space\fi MR }
\providecommand{\MRhref}[2]{%
  \href{http://www.ams.org/mathscinet-getitem?mr=#1}{#2}
}
\providecommand{\href}[2]{#2}

}}
\vspace{.4cm}

%\noindent{\sc Department of Mathematics, University of Copenhagen,
%  Universitets\-par\-ken~5, 2100 Copenhagen {\O}, Denmark}

%\vspace{.2cm}

%\noindent{\sl E-mail address:} {\tt rordam@math.ku.dk}\\
%\noindent{\sl Internet home page:}
%{\tt www.math.ku.dk/$\,\widetilde{\;}$rordam/} \\

%\vspace{.4cm} \noindent \emph{After July 1, 2002:}

\vspace{.2cm}

\noindent{\sc Department of Mathematics, University of Southern
  Denmark, Odense, Campusvej 55, 5230 Odense M, Denmark}

\vspace{.2cm}

\noindent{\sl E-mail address:} {\tt mikael@imada.sdu.dk}\\
\noindent{\sl Internet home page:}
{\tt www.imada.sdu.dk/$\,\widetilde{\;}$mikael/} \\
\end{document}